\documentclass{amsart}
\usepackage{amsfonts}

\theoremstyle{plain}
\newtheorem{acknowledgement}{Acknowledgement}

\numberwithin{equation}{section}
\input{tcilatex}

\begin{document}
\title[Concircular vector fields]{A remarkable property of concircular
vector fields on a Riemannian manifold}
\author[I. Al-Dayel]{Ibrahim Al-Dayel}
\address{Department of Mathematics, College of Science, Imam Muhammad Ibn
Saud Islamic University P.O. Box-65892, Riyadh-11566, Saudi Arabia.}
\email{iaaldayel@imamu.edu.sa}
\author[S. Deshmukh]{Sharief Deshmukh}
\address{Department of Mathematics, College of Science, King Saud
University, P.O.Box-2455, Riyadh-11451, Saudi Arabia}
\email{shariefd@ksu.edu.sa}
\author[O. Belova]{Olga Belova}
\address{Institute of Physical and Mathematical Sciences and Immauel Kant
Baltic Federal University, A. Newsky Str 14, 236216 Kalmingrad, Russia}
\email{obelova@kantiana.ru}
\thanks{This paper is in final form and no version of it will be submitted
for publication elsewhere.}

\begin{abstract}
In this paper, we show that given a nontrivial concircular vector field $%
\boldsymbol{u}$ on a Riemannian manifold $(M,g)$ with potential function $f$%
, there exists a unique smooth function $\rho $ on $M$ that connects $%
\boldsymbol{u}$ to the gradient of potential function \ $\nabla f$, which we
call the connecting function of the concircular vector field $\boldsymbol{u}$%
. Then this connecting function is shown to be a main ingredient in
obtaining characterizations of $n$-sphere $\mathbf{S}^{n}(c)$ and the
Euclidean space $\mathbf{E}^{n}$. We also show that the connecting function
influences topology of the Riemannian manifold.
\end{abstract}

\maketitle

\section{Introduction}

One of important topics in differential geometry of a Riemannian manifold $%
(M,g)$ is studying influence of special vector fields on its geometry as
well as topology. These special vector fields are geodesci vector fields,
Killing vector fields, concircular vector fields, Jacobi-type vector fields
and conformal vector fields on a Riemannian manifold. Moreover, it is well
known that their existence have considerable impact on the geometry of the
Riemannian manifold and these vector fields are used in finding
characterizations of spheres as well as Euclidean spaces (cf. \cite{3}-\cite%
{12}, \cite{15}-\cite{17}, \cite{19}). In \cite{11}, Fialkow initiated the
study of concircular vector fields on a Riemannian manifold. A smooth vector
field $\boldsymbol{u}$ on a Riemannian manifold $(M,g)$ is said to ba a
concircular vector field if%
\begin{equation*}
\nabla _{X}\boldsymbol{u}=fX\text{,\qquad }X\in \mathfrak{X}(M)\text{,}
\end{equation*}%
where $\nabla $ is the Riemannian connection on the Riemannian connection on 
$(M,g)$ and $f:M\rightarrow R$ is a smooth function and $\mathfrak{X}(M)$ is
the Lie algebra of smooth vector fields on $M$ (see also \cite{18}). The
smooth function appearing in the definition of the concircular vector field $%
\boldsymbol{u}$ is called the potential function of the concircular vector
field $\boldsymbol{u}$. A concircular vector field $\boldsymbol{u}$ is said
to be a non-trivial concircular vector field if the potential function $%
f\neq 0$. \ 

Note that a concircular vector field is a closed conformal vector field, a
natural question arises, what is so special about a concircular vector field
among closed conformal vector fields? In this paper, we answer this question
by showing that to each non-trivial concirculalr vector $\boldsymbol{u}$
with potential function $f$ on a connected Riemannian manifold \ $(M,g)$,
there exists a unique smooth function $\rho $ such that $\nabla f=\rho 
\boldsymbol{u}$, where $\nabla f$ is the gradient of the potential function $%
f$. Thus, this unique function $\rho $ connects the gradient $\nabla f$ of
the potential function $f$ and the concircular vector vector field $%
\boldsymbol{u}$ and therefore, we call $\rho $ the \textit{connecting
function} of the concircular vector field $\boldsymbol{u}$. It is
interesting to observe that connecting function $\rho $ is helpful in
finding characterizations of the $n$-sphere $\mathbf{S}^{n}(c)$ as well as
the Euclidean space $\mathbf{E}^{n}$ (cf. theorems, 3.2,-3.4, 4.1).
Moreover, in the last section, we observe that the connecting function $\rho 
$ also influences topology of the Riemannian manifold (cf. theorems 5.1,
5.2).

\section{Preliminaries}

Let $(M,g)$ be an $n$-dimensional Riemannian manifold $(M,g)$ and $%
\boldsymbol{u}$ be a non-trivial concircular vector field on $(M,g)$ with
potential function $f$. \ Then%
\begin{equation}
\nabla _{X}\boldsymbol{u}=fX\text{,\qquad }X\in \mathfrak{X}(M)\text{,} 
\tag{2.1}
\end{equation}%
and the curvature tensor field $R$ of the Riemannian manifold $(M,g)$
satisfies%
\begin{equation}
R(X,Y)\boldsymbol{u}=X(f)Y-Y(f)X\text{,\qquad }X,Y\in \mathfrak{X}(M)\text{,}
\tag{2.2}
\end{equation}%
where%
\begin{equation*}
R(X,Y)Z=\nabla _{X}\nabla _{Y}Z-\nabla _{Y}\nabla _{X}Z-\nabla _{\lbrack
X},_{Y]}Z\text{.}
\end{equation*}%
On an $n$-dimensional Riemannian manifold $(M,g)$, the Ricci tensor $Ric$ is
given by%
\begin{equation}
Ric(X,Y)=\dsum\limits_{i=1}^{n}g\left( R(e_{i},X)Y,e_{i}\right) \text{,} 
\tag{2.3}
\end{equation}%
where $\{e_{1},..,e_{n}\}$ is a local orthonormal frame on $M$. The Ricci
operator $Q$ of the Riemannian manifold $(M,g)$ is a symmetric operator
defined by%
\begin{equation*}
g(QX,Y)=Ric(X,Y)\text{,\qquad }X,Y\in \mathfrak{X}(M)\text{.}
\end{equation*}%
The scalar curvature $S$ of the Riemannian manifold is defined by $S=TrQ$
the trace of the Ricci operator $Q$. The gradient $\nabla S$ of the scalar
curvature satisfies (cf, \cite{1}, \cite{2})%
\begin{equation}
\frac{1}{2}\nabla S=\dsum\limits_{i=1}^{n}\left( \nabla Q\right)
(e_{i},e_{i})\text{,}  \tag{2.4}
\end{equation}%
where the covariant derivative%
\begin{equation*}
\left( \nabla Q\right) (X,Y)=\nabla _{X}QY-Q\nabla _{X}Y\text{.}
\end{equation*}%
Thus, using equations (2.2) and \ (2.3), we conclude%
\begin{equation*}
Ric\left( Y,\boldsymbol{u}\right) =-(n-1)Y(f)\text{.}
\end{equation*}%
Hence, we have%
\begin{equation}
Q(\boldsymbol{u})=-(n-1)\nabla f\text{,}  \tag{2.5}
\end{equation}%
where $\nabla f$ is the gradient of the potential function $f$.

We also have a smooth function $h:M\rightarrow R$ on a Riemannian manifold $%
(M,g)$ associated to concircular vector field $\boldsymbol{u}$, defined by%
\begin{equation}
h=\frac{1}{2}\left\Vert \boldsymbol{u}\right\Vert ^{2}\text{.}  \tag{2.6}
\end{equation}%
Then, using equation (2.1), we find the gradient $\nabla h$ of the smooth
function $h$,%
\begin{equation}
\nabla h=f\boldsymbol{u}\text{.}  \tag{2.7}
\end{equation}%
Note that the Hessian operator $A_{\varphi }$ of a smooth function $\varphi
:M\rightarrow R$ on a Riemannian manifold $(M,g)$, and its Laplacian $\Delta
\varphi $ are defined by%
\begin{equation}
A_{\varphi }X=\nabla _{X}\nabla \varphi \text{,\qquad }\Delta \varphi =\func{%
div}\nabla \varphi =TrA_{\varphi }\text{,}  \tag{2.8}
\end{equation}%
where%
\begin{equation*}
\func{div}X=\dsum\limits_{i=1}^{n}g\left( \nabla _{e_{i}}X,e_{i}\right) 
\text{.}
\end{equation*}%
The Hessian $Hess(\varphi )$ is defined by%
\begin{equation}
Hess(\varphi )(X,Y)=g\left( A_{\varphi }X,Y\right) \text{,\qquad }X,Y\in 
\mathfrak{X}(M)\text{.}  \tag{2.9}
\end{equation}%
Note that if $\varphi $ is a nonconstant smooth function on a compact
Riemannian manifold $(M,g)$ satisfying%
\begin{equation*}
\dint\limits_{M}\varphi =0\text{,}
\end{equation*}%
then the minimal principle, gives%
\begin{equation}
\dint\limits_{M}\left\Vert \nabla \varphi \right\Vert ^{2}\geq \lambda
_{1}\dint\limits_{M}\varphi ^{2}\text{,}  \tag{2.10}
\end{equation}%
where $\lambda _{1}$ is a first nonzero eigenvalue of the Laplace operator $%
\Delta $ acting on smooth functions of $M$.

Recall that the Laplace operator $\Delta $ acting on smooth vector fields on
an $n$-dimensional Riemannian manifold $(M,g)$ is defined by%
\begin{equation*}
\Delta X=\dsum\limits_{i=1}^{n}\left( \nabla _{e_{i}}\nabla _{e_{i}}X-\nabla
_{\nabla _{e_{i}}e_{i}}X\right) \text{,\qquad }X\in \mathfrak{X}(M)\text{,}
\end{equation*}%
where $\{e_{1},..,e_{n}\}$ is an orthonormal frame on $M$. A smooth vector
field $X$ is said to be harmonic if $\Delta X=0$.

\medskip

\bigskip

\section{Connecting functions of concircular vector fields}

In this section, first we show that for a non-trivial conciruclar vector
field $\boldsymbol{u}$ with potential function $f$ on a connected Riemannian
manifold $(M,g)$, there exists a unique smooth function $\rho :M\rightarrow
R $, which we call the connecting function of the conciruclar vector field $%
\boldsymbol{u}$. Then, it is shown that the connecting function $\rho $ can
be used to find characterizations of the $n$-sphere $\mathbf{S}^{n}(c)$ as
well as the Euclidean space $\mathbf{E}^{n}$.

\bigskip

\textbf{Theorem 3.1}: \ \textit{Let }$\boldsymbol{u}$ \textit{be a
non-trivial concircular vector field with potential function }$f$ \textit{on
a connected Riemannian manifold }$(M,g)$. \textit{Then there exists a
uniique function }$\rho :M\rightarrow R$ \textit{satisfying}%
\begin{equation*}
\nabla f=\rho \boldsymbol{u}\text{.}
\end{equation*}

\begin{proof}
Let $\boldsymbol{u}$ be a non-trivial concircular vector field with
potential function $f$ on a connected Riemannian manifold $(M,g)$. Then for
the smooth function $h=\frac{1}{2}\left\Vert \boldsymbol{u}\right\Vert ^{2}$%
, using equation (2.7), and (2.8), we find the following expression for the
Hessian operator $A_{h}$%
\begin{equation*}
A_{h}(X)=X(f)\boldsymbol{u}+f^{2}X\text{.}
\end{equation*}%
Thus, the Hessian $Hess(h)$ of the smooth function $h$ is given by%
\begin{equation}
Hess(h)(X,Y)=X(f)g(\boldsymbol{u},Y)+f^{2}g(X,Y)\text{,\qquad }X,Y\in 
\mathfrak{X}(M)\text{.}  \tag{3.1}
\end{equation}%
Now, as the $Hess(h)$ is symmetric, equation (3.1) implies%
\begin{equation*}
X(f)g(\boldsymbol{u},Y)=Y(f)g(\boldsymbol{u},X)\text{,}
\end{equation*}%
and through which, we conclude%
\begin{equation*}
X(f)\boldsymbol{u}=g(\boldsymbol{u},X)\nabla f\text{.}
\end{equation*}%
Replacing $X$ by $\boldsymbol{u}$ in above equation, we get%
\begin{equation*}
\boldsymbol{u}(f)u=\left\Vert \boldsymbol{u}\right\Vert ^{2}\nabla f\text{,}
\end{equation*}%
which on taking the inner product with $\nabla f$, gives%
\begin{equation}
\left\Vert \boldsymbol{u}\right\Vert ^{2}\left\Vert \nabla f\right\Vert ^{2}=%
\boldsymbol{u}(f)^{2}=g\left( \boldsymbol{u},\nabla f\right) ^{2}\text{.} 
\tag{3.2}
\end{equation}%
Above equation confirms that vector fields $\boldsymbol{u}$ and $\nabla f$
are parallel. Hence, there exists a smooth function $\rho :M\rightarrow R$
such that%
\begin{equation}
\nabla f=\rho \boldsymbol{u}  \tag{3.3}
\end{equation}%
If there is another function $\sigma :M\rightarrow R$, satisfying $\nabla
f=\sigma \boldsymbol{u}$, then we have $\left( \rho -\sigma \right) 
\boldsymbol{u}=0$, which on connected $M$ implies that either $\rho =\sigma $
or $\boldsymbol{u}=0$. However, $\boldsymbol{u}=0$ in equation (2.1), gives $%
f=0$, a contradiction to the fact that $\boldsymbol{u}$ is a non-trivial
concircualr vector field. Hence, $\rho =\sigma $, that is $\rho $ is a
unique function satisfying equation (3.3).
\end{proof}

\bigskip

The unique function $\rho $ guaranteed by Theorem 3.1 that is associated to
the non-trivial concircualar vector field $\boldsymbol{u}$ with potential
function $f$ on a connected Riemannian manifold $(M,g)$ connects the vector
field $\nabla f$ to the vector field $\boldsymbol{u}$. Therefore, we call
the function $\rho $ the \textit{connecting function} of the non-trivial
concircular vector field $\boldsymbol{u}$. In the following results, we show
that the connecting function $\rho $ can be used to characterize a $n$%
-sphere $\mathbf{S}^{n}(c)$ of constant curvature $c$.\bigskip

\textbf{Theorem 3.2}: \textit{An }$n$-\textit{dimensional compact and
connected Riemannian manifold }$(M,g)$ \textit{\ } \textit{admits a
non-trivial concircular vector field }$\boldsymbol{u}$ \textit{with
potential function }$f$ \textit{such that the connecting function }$\rho $ 
\textit{is a constant along the integral curves of }$\boldsymbol{u}$,\textit{%
\ } \textit{if and only if, }$(M,g)$\textit{\ is isometric to the }$n$-%
\textit{sphere }$\mathbf{S}^{n}(c)$.

\begin{proof}
Suppose $(M,g)$ is an $n$-dimensional compact and connected Riemannian
manifold admits a non-trivial concircular vector field $\boldsymbol{u}$ with
potential function $f$ such that the connecting function $\rho $ is a
constant along the integral curves of $\boldsymbol{u}$, that is, $%
\boldsymbol{u}(\rho )=0$. Then, using%
\begin{equation*}
\nabla f=\rho \boldsymbol{u}\text{,}
\end{equation*}%
the Hessian operator $A_{f}$ of the potential function $f$ is computed by
taking covariate derivative in above equation, which is given by%
\begin{equation}
A_{f}X=X(\rho )\boldsymbol{u}+\rho fX\text{,\qquad }X\in \mathfrak{X}(M)%
\text{.}  \tag{3.4}
\end{equation}%
Note that, using equation (3.1), we conclude $\func{div}\boldsymbol{u}=nf$,
\ and integrating this equation leads to%
\begin{equation}
\dint\limits_{M}f=0\text{.}  \tag{3.5}
\end{equation}%
If $f$ is a constant, then above equation concludes that $f=0$, which is
contrary to the assumption that $\boldsymbol{u}$ is a non-trivial circular
vector field. Hence, the potential function $f$ is a nonconstant function.
Now, using the symmetry of the Hessian operator $A_{f}$ in equation (3.4),
we conclude that%
\begin{equation*}
X(\rho )g(\boldsymbol{u},Y)=Y(\rho )g(\boldsymbol{u},X)\text{,}
\end{equation*}%
and it implies that%
\begin{equation*}
X(\rho )\boldsymbol{u}=g(\boldsymbol{u},X)\nabla \rho \text{.}
\end{equation*}%
Replacing $X$ by $\boldsymbol{u}$ in the above equation and using $%
\boldsymbol{u}(\rho )=0$, we conclude%
\begin{equation*}
\left\Vert \boldsymbol{u}\right\Vert ^{2}\nabla \rho =0\text{.}
\end{equation*}%
However, $\boldsymbol{u}\neq 0$ being a non-trivial concircular vector
field, above equation on connected $M$, gives $\nabla \rho =0$, that is, the
connecting function $\rho $ is a constant. Moreover, the constant $\rho $
has to be a nonzero constant, for if $\rho =0$, then Theorem 3.1 will imply $%
f$ is a constant, which is ruled out in the previous paragrph. Taking trace
in equation (3.4), we get $\Delta f=n\rho f$, that is, the nonconstant
function $f$ is eigenfunction of the Laplace operator $\Delta $ acting on
smooth functions on $M$. Since, $M$ is compact, we conclude $n\rho <0$, that
is, the nonzero constant $\rho <0$. We put $\rho =-c$, $c>0$ and we have%
\begin{equation}
\nabla f=-c\boldsymbol{u}\text{ .}  \tag{3.6}
\end{equation}%
Taking covariant derivative in equation (3.6) with respect to $X\in 
\mathfrak{X}(M)$ and using equation (2.1), we get%
\begin{equation}
\nabla _{X}\nabla f=-cfX\text{,\qquad }X\in \mathfrak{X}(M)\text{,} 
\tag{3.7}
\end{equation}%
Hence, the nonconstant function $f$ satisfies the Obata's differential
equation (3.7) (cf. \cite{15}) and thus, the Riemannian manifold $(M,g)$ is
isometric to the sphere $\mathbf{S}^{n}(c)$.

\medskip

Conversely, we know that $\mathbf{S}^{n}(c)$ is a hypersurface of the
Euclidean space $\mathbf{E}^{n+1}$ with unit normal $N$ and the Weingarten
map $A=-\sqrt{c}I$. \ We take a nonzero constant vector field $Z$ on the
Euclidean space $\mathbf{E}^{n+1}$, whose restriction to $\mathbf{S}^{n}(c)$%
, can be expressed as $Z=$ $\boldsymbol{u}+sN$, where $\boldsymbol{u}$ is
the tangential component of $Z$ and $s$ is a smooth function $s=\left\langle
Z,N\right\rangle $ on the sphere $\mathbf{S}^{n}(c)$, and $\left\langle
,\right\rangle $ is the Euclidean metric on $\mathbf{E}^{n+1}$. Taking $X\in 
\mathfrak{X}(\mathbf{S}^{n}(c))$, we get $X(s)=\left\langle Z,\sqrt{c}%
X\right\rangle =\sqrt{c}g(\boldsymbol{u},X)$, where $g$ is the induced
metric on $\mathbf{S}^{n}(c)$. Thus, we conclude%
\begin{equation}
\nabla s=\sqrt{c}\boldsymbol{u}\text{.}  \tag{3.8}
\end{equation}%
Now, as $Z$ is a constant vector field, using the Euclidean connection $D$
on the Euclidean space $\mathbf{E}^{n+1}$, we have $D_{X}Z=0$. For $X,Y\in 
\mathfrak{X}(\mathbf{S}^{n}(c)$, using the Gauss formula for hypersurface $%
D_{X}Y=\nabla _{X}Y-\sqrt{c}g(X,Y)N$, we compute%
\begin{eqnarray*}
0 &=&\left\langle D_{X}Z,Y\right\rangle =X\left\langle Z,Y\right\rangle
-\left\langle Z,\nabla _{X}Y-\sqrt{c}g(X,Y)N\right\rangle \\
&=&Xg(\boldsymbol{u},Y)-g(\boldsymbol{u},\nabla _{X}Y)+\sqrt{c}sg(X,Y) \\
&=&g\left( \nabla _{X}\boldsymbol{u},Y\right) +\sqrt{c}sg(X,Y)\text{,}
\end{eqnarray*}%
and conclude%
\begin{equation}
\nabla _{X}\boldsymbol{u}=-\sqrt{c}sX\text{,\qquad }X\in \mathfrak{X}(%
\mathbf{S}^{n}(c))\text{.}  \tag{3.9}
\end{equation}%
Hence, $\boldsymbol{u}$ is concircular vector field on $\mathbf{S}^{n}(c)$,
with potential function $f=-\sqrt{c}s$, which gives $\nabla f=-\sqrt{c}%
\nabla s$. Using equation (3.8), we conclude%
\begin{equation}
\nabla f=-c\boldsymbol{u}\text{.}  \tag{3.10}
\end{equation}%
Suppose, $f=0$, which will imply $s=0$ and in view of equation (3.8), $%
\boldsymbol{u}=0$, that is, $Z=0$ on $\mathbf{S}^{n}(c)$. As $Z$ is a
constant vector field, we get $Z=0$ on $\mathbf{E}^{n+1}$, which gives a
contradition to the fact that $Z$ is a nonzero constant vector field. Hence, 
$f\neq 0$, that is, $\boldsymbol{u}$ is a non-trivial concircular vector
field with potential function $f$. Then equation (3.10), implies that the
connecting function $\rho =-c$, which is a constant.
\end{proof}

\bigskip

\textbf{Theorem 3.3}: \textit{An }$n$-\textit{dimensional complete and
simply connected Riemannian manifold }$(M,g)$ \textit{\ } \textit{admits a
non-trivial concircular vector field }$\boldsymbol{u}$ \textit{with
potential function }$f$ \textit{such that }$\Delta \boldsymbol{u}=-\lambda 
\boldsymbol{u}$ \textit{for a constant }$\lambda >0$,\textit{\ } \textit{if
and only if, }$(M,g)$\textit{\ is isometric to the }$n$-\textit{sphere }$%
\mathbf{S}^{n}(\lambda )$.

\begin{proof}
Suppose $\boldsymbol{u}$ is a non-trivial concircular vector field on $(M,g)$
with potential function $f$ and connecting function $\rho $ such that $%
\Delta \boldsymbol{u}=-\lambda \boldsymbol{u}$, $\lambda >0$. Using equation
(2.1) and a local orthonormal frame $\{e_{1},..,e_{n}\}$ on $M$, by a
straight forward computation, we get $\Delta \boldsymbol{u}=\nabla f$. Thus,
theorem 3.1, gives $-\lambda \boldsymbol{u}=\rho \boldsymbol{u}$, that is, $%
\left( \rho +\lambda \right) \boldsymbol{u}=0$. Since, a simply connected $M$
is also connected and $\boldsymbol{u}$ being a non-trivial concircular
vector field $\boldsymbol{u}\neq 0$, we must have $\rho =-\lambda $ and
consequently, theorem 3.1, implies $\nabla f=-\lambda \boldsymbol{u}$, which
on using equation (2.1), gives%
\begin{equation}
\nabla _{X}\nabla f=-\lambda fX\text{,\qquad }X\in \mathfrak{X}(M)\text{.} 
\tag{3.11}
\end{equation}%
If potential function $f$ is a constant, then above equation will imply, $%
f=0 $ (as $\lambda >0)$, which is contrary to the assumption that $f$ is
potential function of the non-trivial concircular vector field $\boldsymbol{u%
}$. Hence, equation (3.11) is Obata's differential equation for nonconstant
function $f$ and positive constant $\lambda $, which proves that $(M,g)$ is
isometric to $\mathbf{S}^{n}(\lambda )$.

\medskip

Conversely, on $\mathbf{S}^{n}(c)$, as in the proof of theorem 3.2, there is
a non-trivial concircular vector field $\boldsymbol{u}$ with potential
function $f$ and connecting function $\rho $ that satisfy equations
(3.8)-(3.10), which imply $\Delta \boldsymbol{u}=-c\boldsymbol{u}$.
\end{proof}

\bigskip \medskip

\textbf{Theorem 3.4}: \textit{An }$n$-\textit{dimensional complete and
simply connected Riemannian manifold }$(M,g)$ \textit{\ } \textit{admits a
non-trivial concircular vector field }$\boldsymbol{u}$ \textit{with
potential function }$f$ \textit{and connecting function }$\rho $ \textit{%
satisfying (i) }$g(\nabla f,\nabla \rho )=0$ \textit{and (ii) }$Ric\left(
\nabla f,\nabla f\right) >0$,\textit{\ } \textit{if and only if, }$(M,g)$%
\textit{\ is isometric to the }$n$-\textit{sphere }$\mathbf{S}^{n}(c)$.

\begin{proof}
Suppose $\boldsymbol{u}$ is a non-trivial concicular vector field with
potential function $f$ and connecting function $\rho $ on an $n$-dimensional
Riemannian manifold $(M,g)$ satisfying%
\begin{equation}
g(\nabla f,\nabla \rho )=0\text{\quad and \quad }Ric\left( \nabla f,\nabla
f\right) >0\text{.}  \tag{3.12}
\end{equation}%
Then using theorem 3.1, in the above equations, we conclude%
\begin{equation}
\rho \boldsymbol{u}(\rho )=0\text{\quad and \quad }\rho ^{2}Ric\left( 
\boldsymbol{u},\boldsymbol{u}\right) >0\text{,}  \tag{3.13}
\end{equation}%
that is, $\boldsymbol{u}(\rho )=0$. Using the symmetry of Hessian operator
in equation (3.4), we have%
\begin{equation*}
X(\rho )g(\boldsymbol{u},Y)=Y(\rho )g(\boldsymbol{u},X)\text{,\qquad }X,Y\in 
\mathfrak{X}(M)
\end{equation*}%
and taking $X=$ $\boldsymbol{u}$, in above equation, yields $Y(\rho
)\left\Vert \boldsymbol{u}\right\Vert ^{2}=0$, $Y\in \mathfrak{X}(M)$. As $%
\boldsymbol{u}$ is a non-trivial concircular vector field, we must have $%
Y(\rho )=0$, \ $Y\in \mathfrak{X}(M)$, that is, $\rho $ is a constant and in
view of second equation in equation (3.13), constant $\rho \neq 0$. Now,
equation (2.5) and theorem 3.1, imply%
\begin{equation}
Ric(\boldsymbol{u},\boldsymbol{u})=-(n-1)\rho \left\Vert \boldsymbol{u}%
\right\Vert ^{2}\text{.}  \tag{3.14}
\end{equation}%
Combining equations (3.13) and (3.14), we conclude that the nonzero constant 
$\rho <0$. Taking $\rho =-c$, $c>0$, theorem 3.1, gives $\nabla f=-c%
\boldsymbol{u}$, where $f$ has to be nonconstant, for otherwise we shall
have $\boldsymbol{u}=0$, which is ruled out. Hence, using equation (2.1), we
get the Obata's differential equation%
\begin{equation*}
\nabla _{X}\nabla f=-cfX\text{,\qquad }X\in \mathfrak{X}(M)\text{,}
\end{equation*}%
proving that $(M,g)$ is isometric to $\mathbf{S}^{n}(c)$.

\medskip

Converse trivially follows through the proof of theorem 3.2.
\end{proof}

\section{Characterizations of Euclidean spaces}

In this section, we are interested in finding characterizations of a
Euclidean space using non-trivial concircular vector fields.

\bigskip

\textbf{Theorem 4.1}: \textit{An }$n$-\textit{dimensional complete and
connected Riemannian manifold }$(M,g)$ \textit{admits a non-trivial
concircular vector field }$\boldsymbol{u}$ \textit{with potential function }$%
f$ \textit{\ satisfying }$Ric(\nabla f,\nabla f)=0$ ,\textit{\ if and only
if, }$(M,g)$\textit{\ is isometric to the Euclidean space }$\mathbf{E}^{n}$.

\begin{proof}
Suppose $(M,g)$ is an $n$-dimensional complete and connected Riemannian
manifold $(M,g)$ that admits a non-trivial concircular vector field $%
\boldsymbol{u}$ with potential function $f$, connecting function $\rho $ and
the Ricci curvature satisfies%
\begin{equation}
Ric\left( \nabla f,\nabla f\right) =0\text{.}  \tag{4.1}
\end{equation}%
Using equation (2.5) and theorem 3.1, we have%
\begin{equation*}
Q(\nabla f)=-(n-1)\rho \nabla f\text{,}
\end{equation*}%
which in view of equation (4.1), gives%
\begin{equation}
-(n-1)\rho \left\Vert \nabla f\right\Vert ^{2}=0\text{.}  \tag{4.2}
\end{equation}%
Note that, if $\rho =0$, then theorem 3.1, gives $\nabla f=0$, that is, $f$
is a constant. Thus, as $M$ is connected, equation (4.2), in its both
outcomes, implies that $f$ is a constant. Now, observe that the constant $%
f\neq 0$, owing to the fact that $\boldsymbol{u}$ is non-trivial. Using
equation (2.7), for the function $h=\frac{1}{2}\left\Vert \boldsymbol{u}%
\right\Vert ^{2}$, we find the following expression for its Hessian operator%
\begin{equation*}
A_{h}X=f^{2}X\text{,}
\end{equation*}%
and consequently, we have%
\begin{equation}
Hess(h)=cg\text{,}  \tag{4.3}
\end{equation}%
where $c=f^{2}$ is a nonzero constant. Notice through equation (2.7), that
the function $h$ is not a constant, for if $h$ were to be a constant, as $%
f\neq 0$, it would imply $\boldsymbol{u}=0$ a contradiction. Hence, the
nonconstant function $h$ satisfies equation (4.3) for a nonzero constant $c$%
, proves that the complete and connected Riemannian manifold $(M,g)$ is
isometric to the Euclidean space $\mathbf{E}^{n}$ (cf. Theorem-1, \cite{17}).

Conversely, consider the position vector field%
\begin{equation*}
\boldsymbol{u}=\dsum\limits_{i=1}^{n}x^{i}\frac{\partial }{\partial x^{i}}%
\text{,}
\end{equation*}%
on the Euclidean space $\mathbf{E}^{n}$, where $x^{1},..,x^{n}$ are
Euclidean coordinates. which satisfies $\nabla _{X}\boldsymbol{u}=X$, $X\in 
\mathfrak{X}(\mathbf{E}^{n})$, where $\nabla $ is the Euclidean connection
on $\mathbf{E}^{n}$. Thus, $\boldsymbol{u}$ is a non-trivial concircular
vector field on $\mathbf{E}^{n}$ with potential function $f=1$ and
connecting function $\rho =0$, which satisfies the condition in the
statement of the theorem.
\end{proof}

\bigskip

Our next result shows that harmonic concircular vector fields characterize
Euclidean spaces.

\bigskip

\textbf{Theorem 4.2}: \textit{An }$n$-\textit{dimensional complete and
connected Riemannian manifold }$(M,g)$ \textit{admits a non-trivial
concircular vector field }$\boldsymbol{u}$ \textit{that satisfies }$\Delta 
\boldsymbol{u}=0$,\textit{\ if and only if, }$(M,g)$\textit{\ is isometric
to the Euclidean space }$\mathbf{E}^{n}$.

\begin{proof}
Suppose $\boldsymbol{u}$ is a non-trivial concircular vector field with
potential function $f$ on an $n$-dimensional complete and connected
Riemannian manifold\textit{\ }$(M,g)$, which satisfies $\Delta u=0$. Using
equation (2.1), we compute%
\begin{equation*}
\Delta \boldsymbol{u}=\nabla f\text{.}
\end{equation*}%
Hence, the potential function $f$ is a constant and this constant $f\neq 0$
as $\boldsymbol{u}$ is a non-trivial concircular vector field. Now, equation
(2.7), with $f$ a constant gives%
\begin{equation*}
Hess(h)=cg\text{,}
\end{equation*}%
where $c=f^{2}$ is a nonzero constant. Hence, $(M,g)$ is isometric to the
Euclidean space $\mathbf{E}^{n}$.

Converse is trivial, as the position vector field $\boldsymbol{u}$ on the
Euclidean space $\mathbf{E}^{n}$ is harmonic.
\end{proof}

\section{\protect\bigskip Influence of concircular vector fields on topology}

In this section, we observe that because of the connecting function, we can
exhibit the influence of non-trivial concircular vector fields on topology
of the Reimannian manifolds. Our observations depend on already known
results and therefore results in this sections are simply trivial
applications of known results in differential topology. Recall that by
Reeb's theorem, if a compact smooth manifold $M$ admits a smooth function $%
F:M\rightarrow R$ with exactly two critical points which are non-degenerate,
then $M$ is homeomorphic to an $n$-sphere $S^{n}$. Moreover, it is later
observed by Milnor (cf. theorem 1 p. 166 ,\ \cite{14} ) that this result
holds even if the two critical points are degenerate. Using this modified
Reeb's theorem, we have the following trivial consequence:

\bigskip

\textbf{Theorem 5.1}: \textit{\ If an }$n$-\textit{dimensional compact and
connected Riemannian manifold }$(M,g)$ \textit{admits a non-trivial
concircular vector field }$\boldsymbol{u}$ \textit{\ with potential function 
}$f$ \textit{\ and connecting function }$\rho $ \textit{such that }$\rho
(p)\neq 0$ \textit{for each }$p\in M$\textit{\ and vector field }$%
\boldsymbol{u}$ \textit{\ has only two zeros, then }$M$\textit{\ is
homeomorphic to an }$n$-\textit{sphere}.

\begin{proof}
Using theorem 3.2, we have $\nabla f=\rho \boldsymbol{u}$, and the vector
field $\boldsymbol{u}$ has two zeros say at $p,q\in M$. Then as connecting
function $\rho (x)\neq 0$ on $M$, points $p$, $q$ are critical points of the
potential function \ $f$. Thus, the smooth function $f$ has exactly two
critical points, which proves that $M$ is homeomorphic to $n$-sphere.
\end{proof}

\bigskip

Consider a non-trivial conciruclar vector field $\boldsymbol{u}$ that is
nowhere zero on an $n$-dimensional connected Riemannian manifold $(M,g)$
with potential function $f$ and connecting function $\rho (p)\neq 0$, $p\in
M $. Then by theorem 3.1, the potential function $f$ has no critical points.
If we define a smooth vector field $\xi $ on $M$ by%
\begin{equation*}
\xi =\frac{\nabla f}{\left\Vert \nabla f\right\Vert ^{2}}\text{,}
\end{equation*}%
then, as $\xi (f)=1$, the local flow $\left\{ \phi _{t}\right\} $ of $\xi $
satisfies%
\begin{equation}
f(\phi _{t}(p))=f(p)+t\text{,}  \tag{5.1}
\end{equation}%
which on using escape lemma (cf. \cite{13}), proves that $\xi $ is a
complete vector field and $\left\{ \phi _{t}\right\} $ is the global flow.
Moreover, observe that $f:M\rightarrow R$ is a submersion, consequently, the
lever set $M_{p}=f^{-1}\left\{ f(p)\right\} $ is a compact hypersurface of $%
M $. Now, we have the following:

\bigskip

\textbf{Theorem 5.2}: \textit{\ If an }$n$-\textit{dimensional connected
Riemannian manifold }$(M,g)$ \textit{admits a non-trivial concircular vector
field }$\boldsymbol{u}$, $\boldsymbol{u}(p)\neq 0$, $p\in M$, \textit{\ with
potential function }$f$ \textit{\ and connecting function }$\rho $ \textit{%
such that }$\rho (p)\neq 0$ \textit{for each }$p\in M$\textit{, then }$M$%
\textit{\ is diffeomorphic to }$N\times R$ \textit{for some compact smooth
manifold }$N$.

\begin{proof}
For $p\in M$, we denote by $M_{p}$ the level set $f^{-1}\left\{ f(p)\right\} 
$ of $f$, which is a compact hypersurface of $M$. We define $F:M_{p}\times
R\rightarrow M$ by%
\begin{equation*}
F(q,t)=\phi _{t}(q)\text{,}
\end{equation*}%
which is a smooth map. First, we shall show that $F$ is a surjective: Take $%
x\in M$, then we can find $s\in R$, such that $\phi _{s}(x)=m\in M_{p}$,
with $x=\phi _{-s}(m)$. Consequently,%
\begin{equation*}
F(m,-s)=x\text{.}
\end{equation*}%
Next, we show that $F$ is an injective: Take $(q_{1},t_{1})$, $%
(q_{2},t_{2})\in M_{p}\times R$ such that $F(q_{1},t_{1})=F(q_{2},t_{2})$.
Then we have $\phi _{t_{1}}(q_{1})=\phi _{t_{2}}(q_{2})$ and using equation
(5.1), we get%
\begin{equation*}
f\left( q_{1}\right) +t_{1}=f(q_{2})+t_{2}\text{.}
\end{equation*}%
However, as $q_{1},q_{2}\in M_{p}$, we have $f(q_{1})=f(q_{2})$. Thus, we
get $t_{1}=t_{2}$. and $\phi _{t_{1}}(q_{1})=\phi _{t_{1}}(q_{2})$ implies $%
q_{1}=q_{2}$. Hence $F$ is an injective. Finally, we have%
\begin{equation*}
F^{-1}(x)=(m,-s)=\left( \phi _{s}(x),-s\right)
\end{equation*}%
is also smooth. Hence, $F$ is a diffeomorphism.
\end{proof}

\bigskip \bigskip

\begin{acknowledgement}
This work is supported by King Saud University, Deanship of Scientific
Research, College of Science Research Centre.
\end{acknowledgement}

\end{document}